\documentclass[reqno,10pt]{amsart}
\usepackage{amscd,amssymb,amsmath,color,mathtools}
\usepackage{paralist}
\usepackage{latexsym}
\usepackage[all]{xy}
\usepackage[colorlinks,allcolors=blue]{hyperref}
\usepackage[T1]{fontenc}
\usepackage{color}
\usepackage[english]{babel}
\usepackage{textcomp}
\usepackage{amsmath}
\usepackage{amssymb}
\usepackage{graphicx}
\usepackage{algorithm}
\usepackage{algpseudocode}
\usepackage{hyperref}

\newcommand{\lyxmathsym}[1]{\ifmmode\begingroup\def\b@ld{bold}
  \text{\ifx\math@version\b@ld\bfseries\fi#1}\endgroup\else#1\fi}

\def\..{{,\dots,}}
\def\:{{\colon}}

\def\int{{\rm int}}

\theoremstyle{definition}

\usepackage{booktabs} 

\graphicspath{{figures/}} 
 \usepackage{tikz} 
 \usetikzlibrary{calc}

\usepackage{tikz}
\usetikzlibrary{arrows}

\begin{document}

\author{Or Raz}

\date{\today}

\title{LATTICES OF FLATS FOR SYMPLECTIC MATROIDS}

\address{Einstein Institute of Mathematics, The Hebrew University of Jerusalem, Giv'at Ram, Jerusalem, 91904, Israel}
\email{or.raz1@mail.huji.ac.il}

\begin{abstract}
We are interested in expanding our understanding of symplectic matroids by exploring the properties of a class of symplectic matroids with a "lattice of flats". Taking a well-behaved family of subdivisions of the cross polytope we obtain a construction of lattices, resembling a known definition for the geometric lattice corresponding to ordinary matroid. We construct a correspondence to a set of enveloped symplectic matroids, we denote ranked symplectic matroids. As a by-product of our construction, we also obtain a new way of finding symplectic matroids from ordinary ones and an embedding Theorem into geometric lattices.

The second part of this paper is dedicated to the properties of ranked symplectic matroids and their enveloping ordinary matroids. We focus on establishing a geometric approach to the study of ranked symplectic matroids, demonstrating the ability to take minors, and proving shellability. We finish with a characterization of ranked symplectic matroids using recursive atom orderings.
\end{abstract}

\maketitle

\section{Introduction}

Many books have been written on the subject of matroid theory. We recommend "Matroid Theory" by James Oxley for an introduction to the basics of matroid theory and "Coxeter Matroids" by Borovik, Gelfand, and White for a closer look at symplectic matroids. Developed in order to generalize the concept of independent sets of vector spaces, matroids now have numerous applications in many fields. Some of them include "Information Theory", \cite{[11]}, "Rigidity Theory", \cite{[12]}, and a recent development in algebraic geometry that helped inspire this work, \cite{[13]}. A very useful property of matroids is their large amount of cryptographic characterizations, which contributes a great deal to their many uses. Below is one such characterization using flats:
\subsection{Definition:}
\label{Definition1}

A collection $\mathcal{L}$ of subsets of $\left[n\right]$ is a lattice of flats for a matroid, also called a geometric
lattice if
\begin{enumerate}
\item $\emptyset,\left[n\right]\in\mathcal{L}$
\item $\forall A,B\in\mathcal{L}$ we have $A\cap B\in\mathcal{L}$.
\item For every $A\in\mathcal{L}$ let $\left\{ B_{1},...,B_{m}\right\} $
be the set of elements in $\mathcal{L}$ covering $A$ then $B_{i}\cap B_{j}=A$
for all $i\neq j$, and $\cup_{i=1}^{m}B_{i}=\left[n\right]$.
\end{enumerate}
Some of our readers may be more familiar with the following definition: a poset is a geometric lattice if and only if it is a finite, atomic, graded and submodular lattice; we encourage the readers to use \cite{[1]} for more background on matroid theory.\\

Gelfand and Serganova \cite{[2]} introduced the concept of coxeter matroids, where matroids are a special case in which the coxeter group is the symmetric group. Symplectic matroids are another important example, namely the case where the coxeter group is the hyperoctahedral group, the group of symmetries of the $n$-cube. In contrast to the case of matroids, which will be referred to from now on as ordinary matroids, we know very little about the symplectic case. In this work, we attempt to give a "lattice of flats" characterization for a class of symplectic matroids, which we call ranked. The "lattices of symplectic flats" we construct, denoted $C_n$ lattices, take a similar form to Definition \ref{Definition1} and have many similar properties. Our hope is that we can lay the groundwork for geometric constructions and properties in the spirit of the many simplicial complexes associated with ordinary matroids.\\

Instead of the algebraic property used to define coxeter matroids, we use a more combinatorial approach introduced in \cite{[3]}. We also make use of a characterization of symplectic matroids using independent sets, introduced by T. Chow in \cite{[4]}. Another characterization that exists but is not used in this work is a circuits definition introduced in \cite{[5]} .\\

The structure of this paper is as follows. In Section 2 we define a $C_n$ lattice in a way similar to geometric lattices and introduce some of its basic properties. In Section 3 we present a general discussion of $NBB$ sets which can be of interest and will be the main tool used in Section 4. Section 4 is the heart of this paper that interprets $C_n$ lattices as lattices of flats for ranked symplectic matroids. Another byproduct of this construction is a new way of obtaining symplectic matroids. In section 5 we prove $C_n$ lattices are lexicographically shellable. We end the paper with a characterization of $C_n$ lattices by atom orderings, resembling the way it was done for geometric lattices in \cite{[8]}.\\

Throughout the paper, we use lowercase letters when referring to elements, uppercase letters when referring to sets, and uppercase letters in calligraphic font when referring to collections of sets. For example, we may have $x\in X\in \mathcal{X}$.

\section{Lattice properties}

Geometric lattices were defined on the ground set $[n]$, in the case of $C_n$ lattices, we use a different ground set, the disjoint union of two copies of $[n]$: 

\[
\ensuremath{J=\left[n\right]\cup\left[n\right]^{*}}=\left\{ 1,2...,n\right\} \cup\left\{ 1^{*},2^{*},...,n^{*}\right\} 
\]

We also introduce the map $*:[n] \rightarrow [n]^*$ defined by $i\mapsto i^*$ and $*:[n]^* \rightarrow [n]$ defined by $i^*\mapsto i$. Throughout this paper, we apply $*$ to sets and collections of sets.
\subsection{Definition:}
\label{Definition2}
A set $A\subset P(J)$ is called admissible if $A\cap A^*=\emptyset$ and a maximal admissible set will be called a transversal, in other words, a transversal contains precisely one copy of $i,i^*$ for each $i\in [n]$. We denote the collection of all admissible subsets of $J$ by $P^{ad}(J)$.

\subsection{Definition:}
\label{Definition4}
For a lattice $\mathcal{L}$ and $A,B\in \mathcal{L}$ we say that $B$ covers $A$ and denote $A\lessdot B$  if $A<B$ and for all $C\in \mathcal{L}$ such that $A\leq C\leq B$ we have $C=B$ or $C=A$. $A\in \mathcal{L}$ is an atom if it covers the least element $\hat{0}$.\\

The following is the main subject of this work.

\subsection{Definition:}
\label{Definition3}
A collection $\mathcal{L}$
of subsets of $J$ is a $C_n$ lattice if:
\begin{enumerate}
\item $\emptyset,J\in\mathcal{L}$
\item $\forall A\in\mathcal{L}\setminus \left\{J\right\}$ we have $A\in P^{ad}(J)$.
\item $\forall A,B\in\mathcal{L}$ we have $A\cap B\in\mathcal{L}$.
\item For every $A\in\mathcal{L}$ let $\left\{ B_{1},...,B_{m}\right\} $
be the set of elements in $\mathcal{L}$ covering $A$ then $B_{i}\cap B_{j}=A$
for all $i\neq j$, and $\cup_{i=1}^{m}B_{i}\supseteq J\setminus A^{*}$.
\end{enumerate}
$C_n$ lattices are subsets of the face lattice of the cross-polytope, resembling the ordinary case in which geometric lattices are subsets of the face lattice of the simplex. This gives rise to a neutral realization of $C_n$ lattices.

\subsection*{Example:}
\label{Example1}
The simplest example of a $C_n$ lattice is $\mathcal{L}(U^*_{k,n})$, the collection of admissible sets of size no greater than $k-1$.\\

Before continuing with our main results, we establish some basic properties of $C_n$ lattices. A partially ordered set $S$ in which every pair of elements has a unique supremum (also called join and denoted $\vee$) and a unique infimum (also called meet and denoted $\wedge$) is called a lattice. Every finite lattice contains a greatest element (denoted $\hat{1}$) and a least element (denoted $\hat{0}$).\\

For a $C_n$ lattice we have the following meet and join :

\begin{enumerate}

\item  $A\wedge B=A\cap B$

\item $A\vee B=min\left\{ C\in\mathcal{L}\mid A\cup B\subseteq C\right\} $

\end{enumerate}

As $C_n$ lattices are a finite poset with a meet, they are indeed lattices.

\subsection{Lemma:}
\label{Lemma1}
The atoms in a $C_n$ lattice partition $J$.
\subsubsection*{Proof:}
Using property (1) in \ref{Definition3}, the least element of a geometric lattice is always the empty set. The Lemma now follows directly from property (4) in Definition \ref{Definition3}.

\subsection{Definition:}
\label{Definition5}
A finite lattice is called atomic if every element is a join of atoms.

\subsection{Theorem:}
\label{Theorem1}
$C_n$ lattices are atomic.

\subsubsection*{Proof:}
Let $\mathcal{L}$ be a $C_n$ lattice, we will show that every element $A\in\mathcal{L}$ is
the join of all the atoms contained in it. We have the following:

\[
\vee_{i=1}^{m}A_{i}\subseteq A
\]

As each $A_{i}$ is contained in $A$. We now use Lemma \ref{Lemma1} to
see that every $x\in J$ is in some atom. So, if $x\in A$ we have some atom $B$ s.t. $x\in B$. If $B\setminus A\neq\emptyset$ we will have $\emptyset\subsetneq B\cap A\subsetneq B$ which is a contradiction to $B$ being an atom and so $B\subseteq A$ and we
have:

\[
A\subseteq\vee_{i=1}^{m}A_{i}\Rightarrow A=\vee_{i=1}^{m}A_{i}
\]
\\
We will now show that $C_n$ lattices are graded. In Section $3$ we will prove that there is an underlying family of independent sets which determines this rank function. This should remind our readers of the ordinary case.

\subsection{Definition:}
\label{Definition6}
A lattice $\mathcal{L}$ is called graded if for every $A,B\in \mathcal{L}$ such that $A<B$ and for every two maximal chains $\mathcal{X}=\left\{A=X_1\lessdot ...\lessdot X_m=B\right\}, \mathcal{Y}=\left\{A=Y_1\lessdot ...\lessdot Y_l=B\right\}$ we have $m=l$. This enables us to define the rank function $r(A)$ of a graded lattice to be the length of a maximal chain from $\hat{0}$ to $A$.

\begin{tikzpicture}[x=0.5cm,y=0.5cm]
  \node (max) at (-10,10) {$B$};
  
  \node (b) at (-14.5,4) {$D\neq B$};
  \node (c) at (-10,4) {$D''$};
  \node (d) at (-6,4) {$D'$};

  \node (g) at (-14,1) {$C$};
  \node (h) at (-10,1) {$C''$};
  \node (i) at (-6,1) {$C'$};

  \node (min) at (-10,-1) {$A$};
   
  \draw (min) -- (g) -- (b)
   (min) -- (h) -- (d)
   (min) -- (i) -- (d)
   (min) -- (h) -- (c)
   (min) -- (g) -- (c);
\draw[dashed]  (b) -- (max) 
(d) -- (max)
(c) -- (max);
\end{tikzpicture}\\
    
\subsection{Figure:}
\label{Figure1}
The elements described in the proof of Theorem \ref{Theorem2}. The straight lines represent covering relations, and dashed lines represent containment.

\subsection{Theorem:}
\label{Theorem2}
$C_n$ lattices are graded.

\subsubsection*{Proof:} Let $\mathcal{L}$ be a $C_n$ lattice and $A,B\in\mathcal{L}$ such that $A<B$ . We prove the Theorem by induction on $n$ the
maximum length of a maximal chain between $A$ and $B$.

For length $n=3$, this is obvious as only the empty chain is of length $0$. Let $\mathcal{X},\mathcal{Y}\subseteq \mathcal{L}$
be two maximal chains between $A$ and $B$. 

Let $C\lessdot D$ and $C'\lessdot D'$ be the first two elements other than $A$ in $\mathcal{X},\mathcal{Y}$, respectively, as illustrated in Figure \ref{Figure1}. We first show that $D,D'\neq B$ or $D=D'=B$,  the second case is not interesting as then both chains are of length $3$. Let $D=B$ and assume that there exists $a\in C\setminus C'^*$. By property (4) of Definition \ref{Definition3} we have an element $D''\in \mathcal{L}$, covering $C'$ and containing $a$. As $a\in D''\wedge C$ and $C\lessdot D$, we have
\[
    A\subsetneq D''\wedge C \subset C\Rightarrow C\subset D''\Rightarrow D=D'=D''=B
\]
If there is no such $a$ then let $b\in D'\setminus C'$ If $D'\neq B$ then $b\notin C\cup C^*$. Using property (4) of Definition \ref{Definition3}, we obtain an element $=C''\in \mathcal{L}$ covering $A$ and containing $b$. We have now reduced the problem to the first case, replacing $C'$ with $C''$.\\

To finish the proof of the Theorem, We now first observe the case $(C\setminus A)^*\subset C'$, illustrated in Figure \ref{Figure1}. If that is indeed the case, let $a\in D'\setminus C'$ and define $C''$ as the element that covers $A$ and contains $a$ resulting from property (4) of Definition \ref{Definition3}. As $\mathcal{L}$ is atomic we have $C''<D'$. By applying the induction hypotheses, using the fact $D'\neq B$, to the chains $\left\{A,C',D'\right\}$ and $\left\{A,C'',D'\right\}$, we actually have $C''\lessdot D'$. Once again using property (4) of Definition \ref{Definition3}, we have $C'\cap C''=A$, therefore for each $b\in C\setminus A$ there is an element $D''\in \mathcal{L}$ covering $C''$ and containing $b$. We trivially have
\[
A< D''< B\Rightarrow A \lneq C\wedge D''<C \Rightarrow C<D'' 
\]
We again apply the induction hypotheses, using the fact that $D''\neq B$, to the chains $\left\{A,C,D''\right\}$ and $\left\{A,C'',D''\right\}$ to obtain $C\lessdot D''$.

Let $\mathcal{Z}\subseteq \mathcal{L}$ be some maximal chain from $A$ to $B$ with $C''\lessdot D''$, its first two elements other than $A$.\\

To show $\mathcal{X}$ is of the same length as $\mathcal{Z}$ we apply the induction hypothesis to $\mathcal{X}\setminus \left\{A\right\}$ and $(\mathcal{Z}\setminus \left\{A,C''\right\})\cup \left\{A,C\right\}$. The same is done to show $\mathcal{Y}$ is of the same length as $\mathcal{Z}$, producing the proof.
Finally, We observe that if $(C\setminus A)^*\setminus C'\neq \emptyset$ we can choose $C''=c'$ and $D''$ to be an element covering $C'$ and containing some element in $(C\setminus A)^*\setminus C'$. The rest of the proof is similar.\\

To finish this section we show that $C_n$ lattices are "almost" geometric, with only one missing property, submodularity. 

\subsection{Definition:}
\label{Definition7}
A graded lattice $\mathcal{L}$ will be called submodular if for every $A,B\in \mathcal{L}$ its rank function obeys the identity:
\[
rank(A)+rank(B)\geq rank(A\vee B) +rank(A\wedge B)
\]

\subsection{Theorem:}
\label{Theorem3}
A $C_n$ lattice $\mathcal{L}$ is a geometric lattice if and only if $rank(\mathcal{L})\leq 2$
\subsubsection*{Proof:} 
If $rank(\mathcal{L})\leq 2$, then property (4) of Definition \ref{Definition3} implies property (3) of Definition \ref{Definition1}. As the other properties of Definition \ref{Definition1} are trivially true, we have $\mathcal{L}$ a geometric lattice. If $rank(\mathcal{L})\geq 2$, then there always exist two atoms $A,B$ of $\mathcal{L}$ such that $1\in A$ and $1^*\in B$. We have $A\vee B=J, A\wedge B=\emptyset$ and so:
\[
rank(A)+rank(B)=2\lneq rank(\mathcal{L})=rank(J) +rank(\emptyset)= rank(A\vee B) +rank(A\wedge B)
\]
Therefore, $\mathcal{L}$ is not submodular and is not a geometric lattice.
\section{Independent sets of a lattice}

In \cite{[9]} Andreas Blass and Bruce E. Sagan introduced the concept of $NBB$ sets as a way to calculate the möbius function of a finite lattice. These sets can be thought of as a generalization of $NBC$ sets for ordinary matroids. In this section, we show that if these $NBB$ sets form a matroid, we can order embed the lattice in a geometric lattice, respecting the rank function of the original lattice. We make a small change in the definition of $BB$ sets found in \cite{[9]} by considering only linear orders on the atom set of the lattice, resulting in a simpler definition.

\subsection{Definition:}
\label{Definition8}
Let $(\mathcal{L},\omega)$  be a pair of a finite lattice $\mathcal{L}$ and a linear order $\omega$ on its set of atoms $A(\mathcal{L})$. A nonempty set $D \subseteq A(\mathcal{L})$ of atoms is bounded below or $BB$ if there exist $a\in A(\mathcal{L})$ such that $a$ precedes $D$ in the order $\omega$ and $a\leq \vee D$ in $\mathcal{L}$. A set $B \subseteq A(\mathcal{L})$ is called $NBB$ if $B$ does not contain any bounded below subset.\\

Given a finite lattice $\mathcal{L}$ and $T$ the set of all linear orders on $A(\mathcal{L})$ we define the family of independent sets $I\left(\mathcal{L}\right)$ of $\mathcal{L}$ as follows:

\[
I\left(\mathcal{L}\right)=\left\{ B\in 2^{A(\mathcal{L})}\mid\exists\omega\in T\text{ s.t. }B\text{ is an }NBB\text{ set in }\left(\mathcal{L},\omega\right)\right\}
\]

We refer to the elements in $I\left(\mathcal{L}\right)$ as independent sets of $\mathcal{L}$.  Before continuing, we note the following immediate properties:
\begin{enumerate}
    \item A subset of an independent set is itself independent. In other words, independent sets form an abstract simplicial complex.
    \item A set containing a single atom is not $BB$ in any order.
    \item sets containing at most two atoms are always independent. This is done by choosing any order where one of the two elements is first.
\end{enumerate} 

An immediate implication of the proof of Theorem [1.2] in \cite{[9]} is the following Theorem.

\subsection{Theorem:}
\label{Theorem4}
If $\mathcal{L}$ is a geometric lattice, then $I\left(\mathcal{L}\right)$ is the family of independent sets of the ordinary matroid corresponding to $\mathcal{L}$.\\

For the remainder of this section, we fix a finite, atomic, and graded lattice $\mathcal{L}$.  If $\vee A=X\in \mathcal{L}$ for a set of atoms $A\subseteq A(\mathcal{L})$ we say that $A$ spans $X$ and define $rank(A)\coloneqq rank(X)$.  

\subsection{Lemma:}
\label{Lemma2}
If $I$ is an independent set, then there exists an atom $a\in I$ such that $rank(I\setminus \left\{a\right\})\lneq rank(I)$.

\subsubsection*{Proof:}

For any atom $a\in I$ if $rank(I\setminus \left\{a\right\})= rank(I)$ then $a\leq \vee I=\vee \left(I\setminus\left\{ a\right\} \right)$.

Let $\omega$ be a linear order on $A(\mathcal{L})$ for which $I$ is $NBB$ and let $a$ be the first element of $I$ with respect to $\omega$. we have $a\in \vee (I\setminus \left\{a\right\})$ and so $I\setminus \left\{a\right\}$ is $BB$ by $a$ with respect to $\omega$, which is a contradiction to $I$ being $NBB$ for $\omega$.

\subsection{Lemma:}
\label{Lemma3}
$rank(I)\geq \left|I\right|$ for every independent set $I$.

\subsubsection*{Proof:}

The Theorem follows from Lemma \ref{Lemma2} and an induction argument on the rank of $I$.\\

In the case of geometric lattices, the rank of a flat is the size of a maximal independent set contained in it. The next two lemmas establish the same property for our lattice $\mathcal{L}$.

\subsection{Lemma:}
\label{Lemma4}
Let $I$ be an independent set and $a\in A(\mathcal{L})$ such that $rank(I)\lneq rank(I\cup \left\{a\right\})$ then $I\cup \left\{a\right\}$ is also independent.

\subsubsection*{Proof:}
Let $\omega$ be a linear order on $A(\mathcal{L})$ for which $I$ is $NBB$. Define $\omega_a$ as the linear order obtained from $\omega$ by placing $a$ first. checking the conditions of Definition \ref{Definition8} we see that $I$ is $NBB$ with respect to $\omega_a$ as $a\nleq \vee I$. Therefore, $I\cup \left\{a\right\}$ is $NBB$ with respect to $\omega _a$ as adding the first element does not create new $BB$ subsets.
\subsection{Lemma:}
\label{Lemma5}
The rank of an element $X\in \mathcal{L}$ is the cardinality of a maximal independent set $I$ of atoms smaller than $X$.

\subsubsection*{Proof:}

As all atoms in $I$ are smaller than $X$, we have $\vee I\leq X$, hence Lemma \ref{Lemma3} implies $|I|\leq rank(I)\leq rank(X)$. Next, we find an independent set of $X$ of size $rank(X)$. Let $Y\lessdot X$, as $\mathcal{L}$ is atomic, there exists an atom $a\leq X$ such that $Y\vee a=X$. By an induction argument, we have an independent set $J$ contained in $Y$ with $\left|J\right|=rank(Y)=rank(X)-1$ and as $a\nleq \vee J$ we must have $I=J\cup \left\{a\right\}$ an independent set of size $rank(X)$ contained in $X$.\\

As noted previously, in geometric lattices, the rank of an element is the cardinality of a maximal independent set that it contains. Beginning with the collection of independent sets of a matroid, the elements of its associated geometric lattice - termed flats and ordered by inclusion - are precisely its closed subsets. More formally, a subset $F$ is a flat of the matroid $\mathcal{M}$ (i.e., an element of its geometric lattice) if and only if $rank(F)\lneq rank(F\cup \left\{ e\right\})$ for every element $e\notin F$.

\subsection{Theorem:}
\label{Theorem5}
If $I(\mathcal{L})$ is the family of independent sets of an ordinary matroid, then $\mathcal{L}$ can be uniquely extended to a geometric lattice $\mathcal{P}$, denoted the induced geometric lattice of $\mathcal{L}$, with $I\left(\mathcal{L}\right)=I\left(\mathcal{P}\right)$. In this case, the restriction $f\mid _\mathcal{L}$ of the rank function $f$ of $\mathcal{P}$ is the original rank function on $\mathcal{L}$.

\subsubsection*{Proof:}
We know that there exists a unique ordinary matroid $M$ on the ground set $A(\mathcal{L})$ that corresponds to the family of independent sets $I\left(\mathcal{L}\right)$ of $\mathcal{L}$. We identify each element $X\in \mathcal{L}$ with the set $D=\left\{a\in A(\mathcal{L})\mid a\leq X\right\}$. Using Lemma \ref{Lemma5} we know that $rank_\mathcal{L}(X)$ is the cardinality of a maximal independent set contained in $X$. The same is true for $rank_\mathcal{M}(D)$ by definition. As $I\left(\mathcal{L}\right)=I\left(\mathcal{P}\right)$, we must have $rank_\mathcal{L}(X)=rank_\mathcal{M}(D)$.  

To see $D$ is a flat of $M_\mathcal{L}$ observe that for every atom $a\notin D$ we have $X\lneq \vee (D\cup \left\{a\right\})$ and so:
\[
 rank_\mathcal{L}(D)=rank_\mathcal{L}(X)\lneq rank_\mathcal{L}(D\cup \left\{a\right\})= rank_\mathcal{M}(D\cup \left\{a\right\})
\]

\subsection{Remark:}
\label{Remark1}
It is not always true that $I\left(\mathcal{L}\right)$ is a family of independent sets of an ordinary matroid. A wide class of examples can be found in the study of the adjoint of a matroid. The first example I could find is by A. C. Cheung and can be observed in \cite{[10]}. Another known example is the dual lattice of the Vamos matroid.\\

In addition to the previous embedding Theorem, we can give a characterization of geometric lattices using independent sets. We call an independent set $I$ geometric if $\left|I\right|=rank(I)$.

\subsection{Theorem:}
\label{Theorem6}
$\mathcal{L}$ is a geometric lattice if and only if every independent set $I$ of $\mathcal{L}$ is geometric.

\subsubsection*{Proof:}
If $\mathcal{L}$ is geometric, then $I\left(\mathcal{L}\right)$ is the family of independent sets of a matroid, and so every $I\in I\left(\mathcal{L}\right)$ is geometric. On the other hand, we only need to prove that $\mathcal{L}$ is submodular. Let $X,Y\in \mathcal{L}$ and $I$ be an independent set that spans $X\wedge Y$. Using lemma \ref{Lemma5} we can extend $I$ to independent sets $I\cup J$ that span $X$ and $I\cup J'$ that span $Y$. Note also that $I\cup J\cup J'$ spans $X\vee Y$ and therefore contains an independent set $K$ that spans $X\vee Y$. Using the fact that every independent set is geometric, we obtain the following.
\begin{enumerate}
    \item $rank(X\wedge Y)=\left|I\right|$
    \item $rank(X\vee Y)=\left|K\right|\leq \left|I\right|+\left|J\right|+\left|J'\right|$
\end{enumerate}
This is enough to prove submodularity as follows:
\[
 rank(X)+rank(Y)=2\left|I\right|+\left|J\right|+\left|J'\right|\geq rank(X\wedge Y)+rank(X\vee Y)
\]

In Figure \ref{Figure2} we can see an example of the induced geometric lattice. The partition lattice on the set of $4$ elements is a geometric lattice. Although its dual lattice $(a)$, is not geometric, it is finite, atomic, and graded. Moreover, its family of independent sets constitutes a matroid and therefore can be extended to its induced geometric lattice $(b)$. The element $X$ added to the lattice reduces the rank of the size two subsets of $\left\{ 14/23,13/24,12/34\right\}$, which are the only non geometric independent sets.

\begin{tikzpicture}[x=0.5cm,y=0.5cm]
  \node (max) at (-10,6) {$1/2/3/4$};
  \node (a) at (-17.5,4) {$1/23/4$};
  \node (b) at (-14.5,4) {$14/2/3$};
  \node (c) at (-11.5,4) {$1/24/3$};
  \node (d) at (-8.5,4) {$13/2/4$};
  \node (e) at (-5.5,4) {$12/3/4$};
  \node (f) at (-2.5,4) {$1/2/34$};
  
  \node (g) at (-19,1) {$14/23$};
  \node (h) at (-16,1) {$1/234$};
  \node (i) at (-13,1) {$124/3$};
  \node (j) at (-10,1) {$13/24$};
  \node (k) at (-7,1) {$123/4$};
  \node (l) at (-4,1) {$134/2$};
  \node (m) at (-1,1) {$12/34$};
  \node (min) at (-10,-1) {$1234$};
  
  \node (t) at (-23,2.5) {$(a)$};

  \draw (min) -- (g) -- (a) -- (max) 
   (min) -- (h) -- (c) -- (max)
   (min) -- (i) -- (b) -- (max)
   (min) -- (j) -- (d) -- (max)
   (min) -- (k) -- (e) -- (max)
   (min) -- (l) -- (f) -- (max)
   (min) -- (m) -- (f)
   (g) -- (b)
   (h) -- (a)
   (h) -- (f)
   (i) -- (c)
   (i) -- (e)
   (j) -- (c)
   (k) -- (a)
   (k) -- (d)
   (l) -- (b)
   (l) -- (d)
   (m) -- (e);
   
    \end{tikzpicture}\\

   \begin{tikzpicture}[x=0.5cm,y=0.5cm]
   
    \node (maxx) at (10,6) {$1/2/3/4$};
  \node (aa) at (1,4) {$1/23/4$};
  \node (bb) at (4,4) {$14/2/3$};
  \node (cc) at (7,4) {$1/24/3$};
  \node (dd) at (13,4) {$13/2/4$};
  \node (ee) at (16,4) {$12/3/4$};
  \node (ff) at (19,4) {$1/2/34$};
  \node (x) at (10,4) {$X$};
  
  \node (gg) at (1,1) {$14/23$};
  \node (hh) at (4,1) {$1/234$};
  \node (ii) at (7,1) {$124/3$};
  \node (jj) at (10,1) {$13/24$};
  \node (kk) at (13,1) {$123/4$};
  \node (ll) at (16,1) {$134/2$};
  \node (mm) at (19,1) {$12/34$};
  \node (minn) at (10,-1) {$1234$};
  
 \node (t) at (-3,2.5) {$(b)$};
  
  \draw (minn) -- (gg) -- (aa) -- (maxx) 
   (minn) -- (hh) -- (cc) -- (maxx)
   (minn) -- (ii) -- (bb) -- (maxx)
   (minn) -- (jj) -- (dd) -- (maxx)
   (minn) -- (kk) -- (ee) -- (maxx)
   (minn) -- (ll) -- (ff) -- (maxx)
   (minn) -- (mm) -- (ff)
   (gg) -- (bb)
   (hh) -- (aa)
   (hh) -- (ff)
   (ii) -- (cc)
   (ii) -- (ee)
   (jj) -- (cc)
   (kk) -- (aa)
   (kk) -- (dd)
   (ll) -- (bb)
   (ll) -- (dd)
   (mm) -- (ee)
   (gg) -- (x)-- (maxx)
   (jj) -- (x)
   (mm) -- (x);
   \end{tikzpicture}
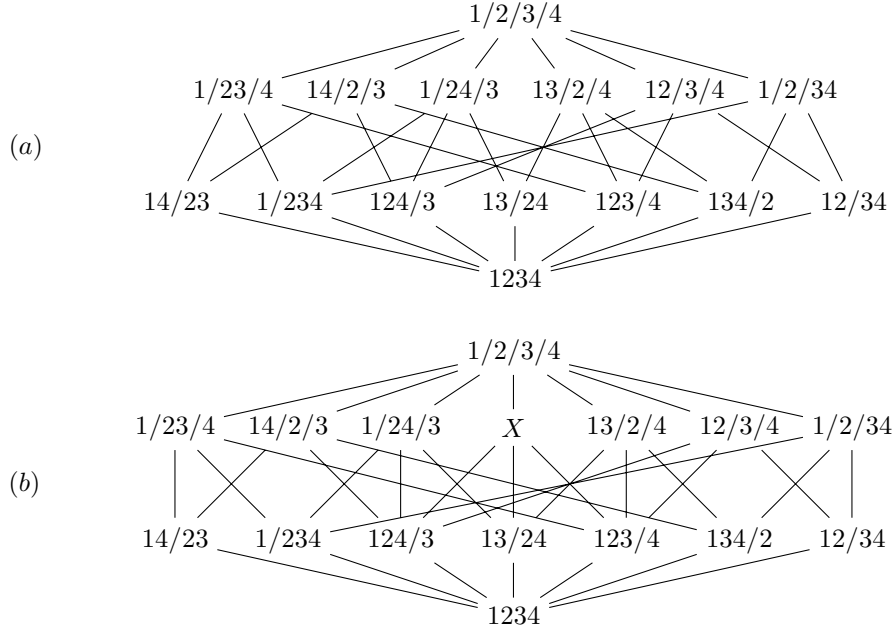
   \subsection{Figure:} \label{Figure2} The dual lattice of the partition lattice on 4 elements $(a)$ and its induced geometric lattice $(b)$.
   
\section{Correspondence to Symplectic matroids}

\subsection{Remark:}
\label{Remark2}
For the purpose of this section we do not consider the full $C_n$ lattice $P^{ad}(J)\cup \left\{ J\right\}$ to be a $C_n$ lattice. The reason is that $P^{ad}(J)\cup \left\{ J\right\}$ is the only $C_n$ lattice in which there are no maximal admissible independent sets. As a result, the symplectic matroid corresponding to $P^{ad}(J)\cup \left\{ J\right\}$ is the same as the one corresponding to the lattice of all admissible subsets of $J$ of size $\leq n-1$. To avoid our correspondence not being injective, we choose to discard $P^{ad}(J)\cup \left\{ J\right\}$.\\

We move on to showing the correspondence between $C_n$ lattices and ranked symplectic matroids, starting with some definitions following \cite{[3]}.\\

An ordering $\leq$ on $J=\left[n\right]\sqcup\left[n\right]^{*}$ is called admissible if and only if $\leq$ is a linear ordering and from $i\leq j$ it follows that $j^*\leq i^*$. A linear order on $J$ induces an order on size $k$ subsets of $J$ in the following way. Let $A=\left\{ a_{1}<a_{2}<...<a_{k}\right\}$ and $B=\left\{ b_{1}<b_{2}<...<b_{k}\right\} $, we set $A\leq B$ if
 \[
 a_{1}\leq b_{1},\,a_{2}\leq b_{2},...,a_{k}\leq b_{k}
 \]

\subsection{Definition:}
\label{Definition9}
A symplectic matroid $\mathcal{B}$ on $J$ is a family of equi-numerous admissible subsets of $J$ called bases, which contains a unique maximal element with respect to every admissible ordering. A symplectic matroid will be called loop-free if every element of $J$ is contained in some basis.\\

An equivalent definition in terms of independent sets was introduced by Timothy Y. Chow in \cite{[4]}:

\subsection{Definition:}
\label{Definition10}
A subset-closed family $\mathcal{I}$ of admissible subsets of $J$ is the family of independent sets of a symplectic matroid if and only if it adheres the following:
\begin{enumerate}
\item For every transversal $T$, $\left\{ I\cap T\mid I\in\mathcal{I}\right\} $ is the family of independent sets of an ordinary matroid with ground set $T$, and
\item  If $I_1$ and $I_2$ are members of $\mathcal{I}$ such that $\left|I_1\right|\lneq\left|I_2\right|$, then either there exists $a\in I_2\setminus I_1$
such that $\left\{ a\right\}\cup I_1 \in \mathcal{I}$ or there exists $a\notin I_1\cup I_2$ such that both $ \left\{ a\right\}\cup I_1 \in \mathcal{I}$ and $(\left\{ a^*\right\}\cup I_1)\setminus {I_2}^* \in \mathcal{I}$.

\end{enumerate}

In this section, we use a slightly different notion of $NBB$ and independent sets defined in section 2. The reason is that we want the independent sets to be subsets of the ground set $J$ and not of the atom set $A(\mathcal{L})$. For a subset $X\subseteq J$ we define $A(X)\subseteq A(\mathcal{L})$ to be $\left\{B\in A(\mathcal{L})\mid B\cap X\neq \emptyset\right\}$, that is, the set of atoms that contain the elements of $X$. We say $X$ is disjoint in $\mathcal{L}$ if the intersection of $X$ with any atom of $\mathcal{L}$ is empty or a singleton. 

Given a linear order $\omega$ on $J$ we also denote by $\omega$ the following linear orders on the set of atoms of a $C_n$ lattice $\mathcal{L}$. $B\lneq _\omega C$ for $B,C \in A(\mathcal{L})$ if $min_\omega B\leq min_\omega C$. These changes will give us the desired connection between the definitions. We now reformulate the results of Section 3 in this setting.

\subsection{Definition:}
\label{Definition11}
Let $(\mathcal{L},\omega)$ be a finite lattice of subsets on the ground set $E$ with a partial order $\omega$ on $E$. A nonempty disjoint in $\mathcal{L}$ set $D \subseteq E$  is bounded below or $BB$ if there exist $a\in E$ such that $a$ is a strict lower bound for all $d\in D$ in the order $\omega$ and $a\in \vee D$.\\
A set $B \subseteq E$ is called $NBB$ if $B$ does not contain any bounded below subset. 

We define the family of independent sets $I\left(\mathcal{L}\right)$ of a finite lattice $\mathcal{L}$ on the ground set $E$ with $T$ the set of all linear orders on $E$.

\[
I\left(\mathcal{L}\right)=\left\{ A\in 2^E\mid\exists\omega\in T\text{ s.t. }A\text{ is an }NBB\text{ set in }\left(\mathcal{L},\omega\right)\right\}
\]

For a $C_n$ lattice $\mathcal{L}$, for which the ground set is $J$, we denote by $I^{ad}\left(\mathcal{L}\right)$ the family of admissible independent sets of $\mathcal{L}$ and by $I^{\overline{ad}}\left(\mathcal{L}\right)$ the family of admissible sets of $\mathcal{L}$ with admissible join.
We are ready to state the main Theorem of this section:

\subsection{Theorem:}
\label{Theorem7}
If $\mathcal{L}$ is a $C_n$ lattice, then $I^{ad}\left(\mathcal{L}\right)$ is the family of independent sets of a symplectic matroid $\mathcal{B}$.\\

\subsection*{Example:}
The Family of independent sets $I(\mathcal{L}(U^*_{k,n}))$ for the $C_n$ lattice in Example \ref{Example1} is the family of all admissible sets of size no grater than $t$.  

To prove this Theorem, we will need to work with two types of ordinary matroids induced by $\mathcal{L}$. We start with the first type of ordinary matroid, corresponding to the family of independent sets of the lattice $\mathcal{L}\cap A=\left\{ X\cap A\mid X\in\mathcal{L}\right\} $ for an admissible $A\in 2^J$. We introduce the following Lemma:
\subsection{Lemma:}
\label{Lemma6}
For any admissible $A\in 2^J$, we have $I\left(\mathcal{L}\right)\cap A = I\left(\mathcal{L}\cap A\right)$.

\subsubsection*{Proof:}
For the inclusion $\subseteq$, it suffices to prove that if $D$ is a minimal $BB$ set in $\mathcal{L}\cap A$ with respect to a partial order $\omega$ on $A$, then $D$ is a $BB$ in $\mathcal{L}$ with respect to any extension of $\omega$ to $J$. This is the case, as the lower bound for $D$ in $\mathcal{L}\cap A$ will also be a lower bound in $\mathcal{L}$. For the inclusion $\supseteq$, let $I$ be $NBB$ in $\mathcal{L}\cap A$ with respect to an order $\omega$, we have that with respect to a partial order starting with $\omega$, $I$ is $NBB$ in $\mathcal{L}$.

\subsection{Corollary:}
\label{corollary1}
For any admissible $A\in 2^J$, we have $I\left(\mathcal{L}\right)\cap A$ the family of independent sets of a matroid.\\

We first describe the inadmissible independent sets in $\mathcal{L}$. For that purpose, we introduce the following Lemma:
\subsection{Lemma:}
\label{Lemma7}
Let $I\in I\left(\mathcal{L}\right)$ and $a\in I$ then exactly one of $I\cup \left\{a^*\right\}\in I\left(\mathcal{L}\right)$ or $rank(I)=rank(\mathcal{L})$ holds.

\subsubsection*{Proof:}
Assume $I\cup \left\{a^*\right\}\notin I\left(\mathcal{L}\right)$, then by Lemma \ref{Lemma4} we have $a^*\in \vee I$. Therefore, $\vee I$ is inadmissible and we get from Definition \
\ref{Definition3} that $rank(I)=rank(\mathcal{L})$.\\

\subsection{Theorem:}
\label{Theorem8}
Let $\mathcal{L}$ be a $C_n$ lattice, then:
\[
I\left(\mathcal{L}\right)=I^{ad}\left(\mathcal{L}\right)\sqcup\left\{ I\cup\left\{ a,a^{*}\right\} \mid I\cup \left\{a\right\},I\cup \left\{a^*\right\}\in I^{\overline{ad}}\left(\mathcal{L}\right) \right\}
\]

\subsubsection*{Proof:}
We first observe that an inadmissible independent set of $\mathcal{L}$ contains exactly one inadmissible pair. If $a,a^*,b,b^*\in B\subseteq A(\mathcal{L})$, then $a\in \vee B\setminus \left\{ a\right\}$ and $b\in \vee B\setminus \left\{ b\right\}$, therefore $B$ is not an independent set.

Next, let $I\cup\left\{ a,a^{*}\right\}\in I\left(\mathcal{L}\right)$ with $I\in I^{ad}\left(\mathcal{L}\right)$. Using Lemma \ref{Lemma7}, we have $rank(I\cup\left\{ a\right\}),rank(I\cup\left\{a^{*}\right\})\neq rank(\mathcal{L})$ and so $I\cup \left\{a\right\},I\cup \left\{a^*\right\}\in I^{\overline{ad}}$. On the other hand, if $I\cup \left\{a\right\},I\cup \left\{a^*\right\}\in I^{\overline{ad}}$, then again by Lemma [4.8] we get $I\cup\left\{ a,a^{*}\right\}\in I\left(\mathcal{L}\right)$.\\

The second type of ordinary matroid that we introduce is the one induced by the original family of independent sets of $\mathcal{L}$, including the inadmissible independent sets. Reformulating Theorem \ref{Theorem5}, we will induce a unique geometric lattice $\mathcal{P}$ on $J$ with the same family of independent sets. For that purpose, we first introduce another useful technical Lemma.

\subsection{Lemma:}
\label{Lemma8}
Let $\mathcal{L}$ be a $C_n$ lattice, if $I\in I^{\overline{ad}}(\mathcal{L})$, then $rank(I)=\left|I\right|$. Moreover, $(I\setminus\left\{a\right\})\cup \left\{a^*\right\}\in I^{\overline{ad}}(\mathcal{L})$ for every $a\in I$.  

\subsubsection*{Proof:}
Let $I\in I^{\overline{ad}}(\mathcal{L})$, if $rank(I)\neq \left|I\right|$ then there exist $I'\subsetneq I$ and $a\in I\setminus I'$ such that $rank(I')=\left|I'\right|$ but $rank(I'\cup \left\{a\right\})\neq\left|I'\cup \left\{a\right\}\right|=rank(I')+1$. Applying Lemma \ref{Lemma3} we must have $rank(I'\cup \left\{a\right\})\gneq rank(I')+1$. As $I'\cup \left\{a\right\}\subseteq I$, we have $I'\cup \left\{a\right\}\in I^{\overline{ad}}(\mathcal{L})$. Therefore, $a,a^*\notin \vee I'$ and by property $4$ of Definition \ref{Definition3} we obtain an element $A\in \mathcal{L}$ covering $\vee I'$ and containing $a$. We must have $\vee (I'\cup \left\{a\right\})=A$ contradicting $rank(I'\cup \left\{a\right\})\gneq rank(I')+1$.

For the second part, let $a\in I\in I^{\overline{ad}}(\mathcal{L})$, we again use property $4$ of Definition \ref{Definition3} to obtain an element $A\in \mathcal{L}$ covering $I\setminus\left\{a\right\}$ and containing $a^*$ with $A\cap (\vee I)=\vee (I\setminus\left\{a\right\})$. We therefore have $\vee((I\setminus\left\{a\right\})\cup \left\{a^*\right\})=A\neq J$ and so $(I\setminus\left\{a\right\})\cup \left\{a^*\right\}\in I^{\overline{ad}}(\mathcal{L})$.
\\

The first step to finding the induced geometric lattice $\mathcal{P}$ is to check the hypothesis of Theorem \ref{Theorem5}, as follows:

\subsection{Theorem:}
\label{Theorem9}
If $\mathcal{L}$ is a $C_n$ lattice, then $I\left(\mathcal{L}\right)$ is a family of independent sets of an ordinary matroid.

\subsubsection*{Proof:}
We prove that $I\left(\mathcal{L},P\right)$ admits the augmentation property for independent sets, stating that for two independent sets $I_1,I_2$ with $\left|I_{1}\right|\lneq\left|I_{2}\right|$ there exists $a\in I_2\setminus I_1$ such that $I_1\cup \left\{a\right\}$ is independent.

Let $I_1,I_2\in I\left(\mathcal{L}\right)$ be such that $\left|I_{1}\right|\lneq\left|I_{2}\right|$. Observe that if $I_1\in I^{\overline{ad}}(\mathcal{L})$, then by Lemmas \ref{Lemma3} and \ref{Lemma8} we have that $rank(I_1)=\left|I_1\right|\lneq \left|I_2\right|\leq rank(I_2)$. Therefore, there exist $a\in I_2$ such that $a\notin \vee I_1$ and by Lemma \ref{Lemma4} $I_1\cup \left\{a\right\}$ is independent.

If $I_1\notin I^{\overline{ad}}(\mathcal{L})$, let $b\in I_1$ be such that $rank(I_1\setminus \left\{b\right\})\lneq rank(I_1)$ using Lemma \ref{Lemma4}. We again use Lemma \ref{Lemma4} to assume that $I_2\in I^{\overline{ad}}(\mathcal{L})$, letting it be smaller but still having $\left|I_{1}\setminus \left\{b\right\}\right|\lneq\left|I_{2}\right|$. We define: 
\[
I_2'=\left\{a\mid (a\in \vee (I_{1}\setminus \left\{b\right\}))\, and \,(\left\{a,a^*\right\}\cap I_2\neq \emptyset)\right\}    
\]
By Theorem \ref{Theorem9} we have $I_2'\in I^{\overline{ad}}(\mathcal{L})$ as it is obtained by a sequence of replacing an element with its star. We have the following:
\[
\left|I_{2}'\right|=rank(I_{2}')\leq rank(I_{1}\setminus \left\{b\right\})=\left|I_{1}\setminus \left\{b\right\}\right|\lneq \left|I_{2}\right|
\]
Therefore, there exists $a\in I_2$ such that $a,a^*\notin \vee (I_{1}\setminus \left\{b\right\})$. As $rank(I_{1}\setminus \left\{b\right\})\lneq rank(I_{2})\lneq rank(\mathcal{L})$ and $\vee ((I_{1}\setminus \left\{b\right\})\cup \left\{a\right\})$ covers $\vee(I_{1}\setminus \left\{b\right\})$ we have $(I_{1}\setminus \left\{b\right\})\cup \left\{a\right\}\in I^{\overline{ad}}(\mathcal{L})$. Finally:
\[
rank(((I_{1}\setminus \left\{b\right\})\cup \left\{a\right\}))\lneq rank(\mathcal{L})=rank(I_1)\leq rank(I_1\cup \left\{a\right\})
\]
then by Lemma \ref{Lemma4} $I_1\cup \left\{a\right\}\in I(\mathcal{L})$.\\

\subsubsection*{Proof of Theorem \ref{Theorem7}:}
By definition, $NBB$ is a hereditary property, so $I\left(\mathcal{L}\right)$ is a subset closed family. 

\begin{enumerate}

\item Given a transversal $T\subset J$ we have by Lemma \ref{Lemma6} that $I\left(\mathcal{L}\right)\cap T=I\left(\mathcal{L}\cap T \right)$. By Theorem \ref{Theorem5} we then find that $I\left(\mathcal{L},P\right)\cap T$ is the family of independent sets of an ordinary matroid on the ground set $T$.

\item Let $I_1,I_2\in I^{ad}\left(\mathcal{L}\right)$ be such that $\left|I_1\right|\lneq \left|I_2\right|$. Assume that for every $a\in I_2\setminus I_1$ such that $ I_1\cup \left\{ a\right\}\in I\left(\mathcal{L}\right)$ has $a^*\in I_1$. By Theorem \ref{Theorem9}, there exists such an $a$, and so by Lemma \ref{Lemma7} we have $I_1\in I^{\overline{ad}}\left(\mathcal{L}\right)$. We now observe that $\left|I_1\setminus I_2^*\right|\lneq \left|I_2\setminus I_1^*\right|$ for the smaller independent sets $I_1\setminus I_2^*,I_2\setminus I_1^*$. By Theorem \ref{Theorem9} we have some $b\in I_2\setminus I_1^*$, $b\notin I_1\setminus I_2^*$ such that $(I_1\setminus I_2^*)\cup \left\{ b\right\}\in I\left(\mathcal{L}\right)$. We actually have $(I_1\setminus I_2^*)\cup \left\{ b\right\}\in I^{ad}\left(\mathcal{L}\right)$ because we have removed all inadmissible pairs in $I_1\cup I_2$. We have found $b\in I_2\setminus I_1$ with $b^*\notin I_1\cup I_2$. By our initial assumption $I_1\cup \left\{ b\right\}\notin I\left(\mathcal{L}\right)$, and as $I_1\in I^{\overline{ad}}\left(\mathcal{L}\right)$ we have $I_1\cup \left\{ b^*\right\}\in I\left(\mathcal{L}\right)$ by Lemma \ref{Lemma4}. As $b\notin I_1$, we get $I_1\cup \left\{ b^*\right\}\in I^{ad}\left(\mathcal{L}\right)$ satisfying the second part of condition (2) of Definition \ref{Definition10}.

\end{enumerate}

After showing that every $C_n$ lattice corresponds to a symplectic matroid, we continue by identifying the symplectic matroids corresponding to $C_n$ lattices.
\subsection{Theorem:}
\label{Theorem10}
Let $\mathcal{L}$ be a $C_n$ lattice, and $\mathcal{P}$  denote the geometric lattice induced by $\mathcal{L}$. Then $\mathcal{P}$ contains no additional admissible elements beyond those already in $\mathcal{L}$; that is, if $X \in \mathcal{P}$ is admissible, then necessarily $X \in \mathcal{L}$.

\subsubsection*{Proof:}
We first prove that if $I\in I^{ad}(\mathcal{L})\setminus I^{\overline{ad}}(\mathcal{L})$ then $\vee_\mathcal{P}I$ is inadmissible. 

Using Lemma \ref{Lemma2}, let $a\in I$ be such that $I\setminus \left\{a\right\}\in I^{\overline{ad}}(\mathcal{L})$. As a result of property 4 of Definition \ref{Definition3} we must have $a^*\in \vee_\mathcal{L}(I\setminus \left\{a\right\})$. Using Lemma \ref{Lemma8} we have
\[
rank_\mathcal{L}(I\setminus \left\{a\right\})=\left|I\setminus \left\{a\right\}\right|=rank_\mathcal{P}(I\setminus \left\{a\right\})\Rightarrow a^*\in \vee_\mathcal{P}(I\setminus \left\{a\right\})\subseteq X
\]

To prove the Theorem, let $X$ be an admissible element of $\mathcal{P}$, and let $I$ be an independent set such that $\vee_\mathcal{P}I=X$. We must then have $I\in I^{\overline{ad}}(\mathcal{L})$ and by Lemma \ref{Lemma8}, $rank_\mathcal{L}(I)=\left|I\right|=rank_\mathcal{P}(I)$. As $\mathcal{P}$ is closed under intersection, we have 
\[
I\subseteq(\vee_\mathcal{P}I)\cap(\vee_\mathcal{L}I)\in \mathcal{P} \Rightarrow  rank_\mathcal{P}((\vee_\mathcal{P}I)\cap(\vee_\mathcal{L}I))=rank_\mathcal{P}(X)\Rightarrow X=\vee_\mathcal{L}I\in \mathcal{L} 
\]
\\
We can now use Theorems \ref{Theorem8}, \ref{Theorem9}, and \ref{Theorem10} to formulate the characterization of symplectic matroids corresponding to $C_n$ lattices using rank functions. 

\subsection{Definition:}
\label{Definition12}
Let $\mathcal{I}$ be the family of independent sets of a symplectic matroid. We call $I\in \mathcal{I}$ a delta-independent set if $I\cup\left\{ a\right\}\in \mathcal{I}$ or $I\cup\left\{ a^{*}\right\}\in \mathcal{I}$ for every $a\in J$. We denote the family of delta-independent sets by $\mathcal{I}^{\triangle}$.

\subsection{Theorem:}
\label{Theorem11}
Let $\mathcal{I}$ be a family of independent sets of a symplectic matroid of rank $d$ with no loops, then $\mathcal{I}$ corresponds to a $C_n$ lattice and denoted a ranked symplectic matroid if and only if the function $r:P\left(J\right)\rightarrow\mathbb{Z}_{+}$ is the rank function of a loop-free ordinary matroid on the ground set $J$ for the following function $r$:

\[
r\left(A\right)=min\left(d,\max_{I\subseteq A}\begin{cases}
\left|I\right|+2 & I\in\mathcal{I}^{\triangle}\land\exists\left\{ a,a^{*}\right\} \subseteq A\setminus I,I\cup\left\{ a\right\} ,I\cup\left\{ a^{*}\right\} \in\mathcal{I}\\
\left|I\right| & I\in\mathcal{I}
\end{cases}\right)
\]

\subsubsection*{Proof:}
The "only if" direction is exactly Theorem \ref{Theorem8} and \ref{Theorem9}, as the rank of a set in an ordinary matroid is the size of a maximal independent set contained in it. For the "if" direction it is enough to show, again using Theorem \ref{Theorem9}, that the admissible part of $\mathcal{P}$ is a $C_n$ lattice for the induced geometric lattice $\mathcal{P}$ of $r$. We note that every atom of $\mathcal{P}$ is admissible by the definition of $r$ and $\mathcal{I}$ being without loops. Therefore, we are reduced to showing that any non-maximal admissible flat $F$ of $\mathcal{P}$ is covered by exactly one inadmissible flat $F\cup F^*$. To see that any inadmissible flat covering $F$ is contained in $F\cup F^*$, let $I\cup	\left\{ a\right\}$ be an independent set spanning $F$ such that $I\cup	\left\{ a^*\right\}$ is also independent. If $b,b^*\notin F$, then we must have $b,b^*\notin \vee_r(I\cup	\left\{ a,a^*\right\})$, with $\vee_r$ being the join in the geometric lattice corresponding to $r$. As both $I\cup \left\{ a,b\right\}$ and $I\cup \left\{ a,b^*\right\}$ are independent and $I\cup \left\{ a\right\}\in \mathcal{I}^{\triangle}$, contradicting $\vee_r(I\cup	\left\{ a,a^*\right\}),\vee_r(I\cup	\left\{ b,b^*\right\})$ being of the same rank. As the rank function $r$ does not allow independent sets with more than one pair of inadmissible elements, we also must have $F\cup F^*$ contained in the inadmissible flat covering $F$. \\ 

Using ordinary matroid terminology, we have defined a "cryptomorphism" $f$ sending a $C_n$ lattice to the symplectic matroid defined by its family of admissible independent sets. By Theorems \ref{Theorem5} and \ref{Theorem8} the family of admissible independent sets is uniquely determined by the induced geometric lattice, and so we have $f$ injective; see Remark \ref{Remark2}.

The inverse of $f$ is obtained by taking the sublattice of admissible flats of the geometric lattice corresponding to the rank function of the ranked symplectic matroid.

The next two Remarks now easily follow:

\subsection{Remark 1:}
\label{Remark3}
An ordinary matroid on $J$ of rank $\geq 3$ with admissible bases, and so also a symplectic matroid, is never a ranked symplectic matroid. The rank function from Theorem \ref{Theorem11} will never be submodular.
\subsection{Remark 2:}
 \label{Remark4}
In Theorem \ref{Theorem11} we have required $\mathcal{I}$ to be a family of independent sets of a symplectic matroid, we can weaken this condition, the same proof will work for the following Theorem:

Let $\mathcal{M}$ be an ordinary matroid of rank $d$ on the ground set $J$. If $\mathcal{M}$ contains an admissible base and its rank function is the one described in Theorem \ref{Theorem11} then $\mathcal{L}$, the admissible part of the geometric lattice $\mathcal{P}$ corresponding to $\mathcal{M}$, is a $C_n$ lattice. Moreover, the ranked symplectic matroid corresponding to $\mathcal{L}$ is the family of admissible bases of $\mathcal{M}$. \\

Remark \ref{Remark4} introduces a new way of constructing families of symplectic matroids using ordinary matroids. This can be helpful because the theory of ordinary matroids is much more developed. We show an example of how any spike with no tip gives rise to a ranked symplectic matroid.

In \cite{[6]} T. Zaslavsky introduced lift matroids on bias graphs for which spikes are a special case.

\subsection{Definition:}
\label{Definition13}
A biased graph is a pair $(G,\mathcal{C})$ where $G$ is a finite undirected multigraph and $\mathcal{C}$ is a set of cycles of $G$ satisfying the following theta property.

\subsubsection*{Theta property:} For every two cycles $C_1$ and $C_2$ in $\mathcal{C}$ that intersect in a nonempty path, the third cycle in $C_1\cup C_2$ is also in $\mathcal{C}$.\\

The cycles of $\mathcal{C}$ are called balanced and the other cycles are unbalanced. The lift matroid represented by $\mathcal{M}(G,\mathcal{C})$ has for its ground set the set of edges of $G$ and, as independent sets, the sets of edges containing at most one cycle, which is unbalanced. A spike with no tip is a lift matroid with $G$ the graph obtained from the cycle on $n$ vertices by doubling every edge and $\mathcal{C}$ some set of cycles that do not contain any pair of double edges and satisfying the theta property. \\

\subsection{Theorem:}
\label{Theorem12}
The family of admissible bases of any spike with no tip, $\mathcal{M}(G,\mathcal{C})$, constitutes a ranked symplectic matroid.

\subsubsection*{Proof:}

We will show that $\mathcal{M}(G,\mathcal{C})$ is a simple matroid admitting the properties described in Remark \ref{Remark4}.

First, notice that $J=[n]\sqcup [n^*]$ with $n$ being the size of the original cycle in $G$, labeling one of the edges $\left\{x,y\right\}$ with $i\in [n]$ gives the other $\left\{x,y\right\}$ edge the label $-i$.

Every simple cycle in $G$ is a set of the form $\left\{ i,i^*\right\}$ for some edge $i$ or a transversal $T$. Therefore, every maximal independent set of $\mathcal{M}(G,\mathcal{C})$ is an unbalanced transversal $T$ or a missing transversal of the form $(T\setminus \left\{i\right\})\cup \left\{j\right\}$, for a transversal $T$, and $i,j^*\in T$.

To see that there always exists an unbalanced transversal, we observe that, for example, $[n],([n]\setminus \left\{1\right\})\cup \left\{1^*\right\}$ cannot be both balanced. Otherwise, we will have $\left\{1,1^*\right\}$ balanced by the Theta property. Thus, we have found an admissible basis for $\mathcal{M}(G,\mathcal{C})$.

We are left to show that the rank function of $\mathcal{M}(G,\mathcal{C})$ agrees with the rank function $r$ defined in Theorem \ref{Theorem11}. Recall that for an ordinary matroid, the rank of a set is the size of a maximal independent set contained in it. If $X\subseteq J$ is admissible, then trivially $rank_{\mathcal{M}(G,\mathcal{C})}(X)=r(X)$. This is also the case if $X$ contains an unbalanced transversal. Otherwise, let $I\subseteq X$ be a maximal admissible independent set, if $I^*\cap X\neq \emptyset$ then for $i\in I^*\cap X$  we have $I\cup \left\{i\right\}$ independent in $\mathcal{M}(G,\mathcal{C})$ and $rank_{\mathcal{M}(G,\mathcal{C})}(X)=\left|I\right|+1$. Every non-maximal admissible independent set of $\mathcal{M}(G,\mathcal{C})$ is delta-independent we have:
\[
rank_{\mathcal{M}(G,\mathcal{C})}(X)=\left|I\right|+1=\left|I\setminus \left\{i\right\}\right|+2=r(X)
\]
\\

\begin{tikzpicture}[x=0.5cm,y=0.5cm]
  \node (a) at (-13,0) {$v$};
  \node (b) at (-13,6) {$u$};
   \node (c) at (-16,3) {$1$};
  \node (d) at (-14,3) {$1^*$};
   \node (e) at (-10,3) {$2$};
  \node (f) at (-12,3) {$2^*$};
  \draw [red] (a) to[out=0,in=0] (b);
   \draw [blue]  (b) to[out=180,in=180] (a);
     \draw [blue] (a) to[out=40,in=-40] (b);
   \draw [red] (b) to[out=220,in=140] (a);
   
  \node (max) at (-5,4) {$\{1,2,1^*,2^*$\}};
  \node (f') at (-7,2) {$\{1,2^*\}$};
  \node (e') at (-3,2) {$\{1^*,2\}$};
  \node (min) at (-5,0) {$\emptyset$};
  
  \draw (min) -- (f') -- (max) 
   (min) -- (e') -- (max);
   
   \draw (2, 0) rectangle (6, 4);
 \node (k) at (7,4.5) {$\{1,2\}$};
 \node (k') at (1,-0.5) {$\{1^*,2^*\}$};
 \draw [blue] (2, 0) -- (6, 4);
 
   \node (a) at (-13,-2) {($a$)};
  \node (b) at (-5,-2) {($b$)};
  \node (c) at (4,-2) {($c$)};
   \end{tikzpicture}
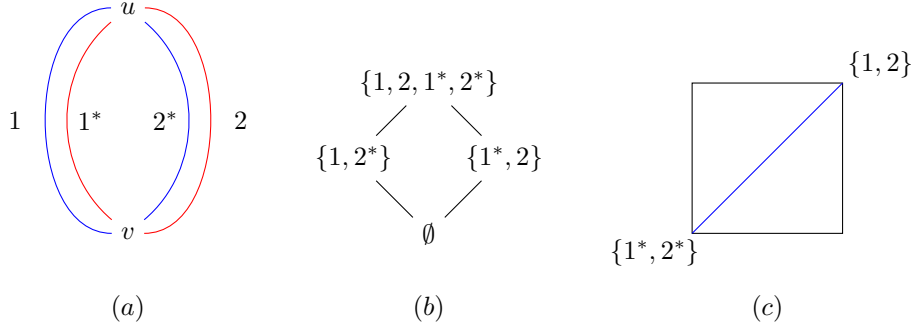
   \subsection{Figure:}
   \label{Figure3}
   The graph defining a spike with no tip, with the blue and red cycles being balanced $(a)$, the $C_n$ lattice induced by it $(b)$, and the convex polytopes associated with its ranked symplectic matroid $(c)$.\\
   
   In \cite{[4]} T. Chow introduced graphic symplectic matroids. We note that our construction of spikes with no tip is not necessarily a graphic symplectic matroids, as the transversal being even does not correspond to the number of starred elements.

\section{Shellability}
The purpose of this section is to prove that the $C_n$ lattices are shellable. We start with a few definitions. 

\subsection{Definition:}
\label{Definition14}
 A simplicial complex $\Delta$ is a set of simplices that satisfy the following conditions:
 \begin{enumerate}
 \item Every face of a simplex from $\Delta$ is also in $\Delta$.
 \item The nonempty intersection of any two simplices $\sigma,\tau\in \Delta$ is a face of both $\sigma$ and $\tau$.
 \end{enumerate}
 
\subsection{Definition:}
\label{Definition15}
Let $\Delta$ be a finite simplicial complex. We say that $\Delta$
is pure $d$ dimensional if all its facets (inclusion maximal faces) are of dimension $d$.
A pure $d$ dimensional simplicial complex $\Delta$ is said to
be shellable if its facets can be ordered $F_{1},...,F_{t}$ in such
a way that:
\[
\overline{F_{j}}\cap\cup_{i=1}^{j-1}\overline{F_{i}}
\]
 Is a pure $\left(d-1\right)$ dimensional complex for $j=2,3,...,t$,
with: 
\[
\overline{F_{j}}=\left\{ G\,\text{a simplicial complex}\mid G\subseteq F_{j}\right\} 
\]
. 

To a finite poset $\mathcal{L}$ one can associate a simplicial complex $\Delta\left(\mathcal{L}\right)$
(the order complex of $\mathcal{L}$) of all chains of $\mathcal{L}$. Also, if $\mathcal{L}$
is a graded poset of rank $d$, then $\Delta\left(\mathcal{L}\right)$ is pure
$d$ dimensional. We say that a finite and graded poset $\mathcal{L}$ is shellable
if its order complex $\Delta\left(\mathcal{L}\right)$ is shellable.

We now define a recursive atom ordering.

\subsection{Definition:}
\label{Definition16}
A graded poset $\mathcal{L}$ is said to admit a recursive atom ordering if
the rank of $\mathcal{L}$ is 1 or if the rank of $\mathcal{L}$ is greater than $1$
and there is an ordering $A_{1},...,A_{t}$ of the atoms of $\mathcal{L}$ that
satisfies:
\begin{enumerate}
\item For all $j=2,...,t$ we have $\left[A_{j},\hat{1}\right]$ admit a
recursive atom ordering in which the atoms of $\left[A_{j},\hat{1}\right]$
that come first in the ordering are those that cover some $A_{i}$
where $i<j$.
\item For all $i<j,$ if $A_{i},A_{j}<B$ for some $B\in \mathcal{L}$ then there
is some $k<j$ and an element $B\geq C\in \mathcal{L}$ such that $C$ covers $A_{k}$
and $A_{j}$. 
\end{enumerate}

In \cite{[7]} A. Bjorner and M. Wachs proved the following Theorem:

\subsection{Theorem:}
\label{Theorem13}
Let $P$ be a graded poset. If $P$ admits a recursive atom ordering, then the order complex $\Delta\left(P\right)$ is shellable.\\

We are now ready to start proving that $C_n$ lattices are shellable, using induction on the number of atoms. We start by introducing some technical Lemmas.

\subsection{Remark:}
\label{Remark5}
It is possible to prove a stronger type of lexicographic shellability called EL-shellability by taking the known EL-labeling for geometric lattices (can be found in \cite{[10]}) of its enveloping ordinary matroid, starting with the atoms of an admissible coatom. To expand the range of shelling orders and support the characterization in part 6, we chose to establish CL-shellability, which is equivalent to admitting a recursive atom ordering.

\subsection{Lemma:}
\label{Lemma9}
If $A$ is an atom of a $C_n$ lattice $\mathcal{L}$ of rank $\geq 3$ then $A^{*}$ is also an atom of
$\mathcal{L}$.

\subsubsection*{Proof:}

Assume that $A,B$ are atoms of $\mathcal{L}$ with $B^{*}\cap A\neq\emptyset$.
Assume by negation that $A\setminus B^*\neq\emptyset$ and let $a\in A\setminus B^{*}$. Applying property $4$ of Definition \ref{Definition3} we have $a\notin B$ and an element $C\in \mathcal{L}$ containing $a$ and covering $B$. As $a\in C\cap A\in \mathcal{L}$, we have $C\cap A=A\Rightarrow A\subseteq C$. By our assumption, we now have:
\[
A\cap B^{*}\neq\emptyset\overset{B\subseteq C}{\Rightarrow}A\cap C^{*}\neq\emptyset\overset{A\subseteq C}{\Rightarrow}C\cap C^{*}\neq\emptyset
\]
That is, $C$ is inadmissible. As $rank(C)=2\lneq rank(\mathcal{L)}$ and $J$ is the only inadmissible element in $\mathcal{L}$, this is a contradiction. Therefore, $B^{*}\cap A\neq\emptyset\Rightarrow B^{*}=A$.

\subsection{Lemma :}
\label{Lemma10}
Let $\mathcal{L}$ be a $C_n$ lattice and $A\in \mathcal{L}$ then the restriction $\left[A,J\right]\subseteq\mathcal{L}$ is order isomorphic
to a $C_n$ lattice on the ground set $J'=J\setminus\left(A\cup A^{*}\right)$.

\subsubsection*{Proof:}

Let $A\in \mathcal{L}$, we will work with the lattice morphism:

\[
\psi:\left[A,J\right]\rightarrow 2^{J'}
\]

\[
X\mapsto\begin{cases}
J' & X=J\\
X\setminus A & X\neq J
\end{cases}
\]

We observe that $\psi(X)\subseteq J'$ for any $X\in \left[A,J\right]$ because $X\cap A^*=\emptyset$ or $X=J$. It is also clear that $\psi$ is injective, and hence it is a lattice isomorphism on its image.

We continue by proving the four properties of Definition
\label{Definition29}:
\begin{enumerate}
\item $\emptyset,J'\in\mathcal{L}$ as we have $A\mapsto \emptyset$
and $J\mapsto J'$.
\item Every element $B\in \psi(\left[A,J\right])$ except for $J'$ is
admissible as it is a subset of a proper element of $\mathcal{L}$.
\item $\forall B,C\in\psi(\left[A,J\right])$ we have: 
\[
A\subseteq\left(A\cup B\right)\cap\left(A\cup C\right)\in\left[A,J\right]\Rightarrow B\cap C\in\psi(\left[A,J\right])
\]
\item For $J'\neq B\in\psi(\left[A,J\right])$ let $\left\{ B_{1},...,B_{m}\right\}$ be
the set of elements in $\psi(\left[A,J\right])$ covering
$B$. If $\left\{ B_{1},...,B_{m}\right\}=\emptyset$ or $\left\{ B_{1},...,B_{m}\right\}=J'$, then property $(4)$ follows directly. We now consider the nontrivial case where $|\left\{ B_{1},...,B_{m}\right\}|\geq 2$. We know that $\psi$ is a lattice isomorphism from an upper interval, hence for any $X\in \left[A,J\right]$, we have $Y$ cover $X$ in $\mathcal{L}$ if and only if $\psi(Y)$ covers $\psi(X)$ in $\psi(\left[A,J\right])$.
Therefore, we have $\psi ^{-1}(B)=B\cup A\in\mathcal{L}$ and a set
of elements $\left\{ \psi(B_{1}),...,\psi(B_{m})\right\}= \left\{ B_{1}\cup A,...,B_{m}\cup A\right\}$ covering $B\cup A$
in $\mathcal{L}$.
As a result of property $(4)$ for $\left\{ B_{1}\cup A,...,B_{m}\cup A\right\}$ in $\mathcal{L}$ we have:
\[
B_{i}\cap B_{j}=((B_{i}\cup A)\cap (B_{j}\cup A))\setminus A=(B\cup A)\setminus A=B
\]
And also:
\[
\cup_{i=1}^{m}B_{i}=(\cup_{i=1}^{m}(B_{i}\cup A))\setminus A=(J\setminus A^*)\setminus A=J'
\]
\end{enumerate}

\subsection{Lemma:}
\label{Lemma11}
If the number of atoms in a $C_n$ lattice $\mathcal{L}$ is $t$, then the number of atoms of $\left[A,J\right]\subseteq\mathcal{L}$
for $A$ an atom of $\mathcal{L}$ is $\leq t-2$.

\subsubsection*{Proof:}

As $\mathcal{L}$ is atomic, we have that each of its elements is a join of some of its atoms. By property (4) of $C_n$ lattices the elements that cover $A$ pairwise disjoint (excluding the elements of $A$) and so every element covering $A$ is a choice of some atoms
from a set of size $t-2$, namely the set $A(\mathcal{L})\setminus \left\{A,A^*\right\}$. There can be at most $t-2$ choices.\\
 
To finish our technical preparations, we observe that we can always work with $C_n$ lattices in which every atom is a singleton. This is done for geometric lattices by deleting all the elements in an atom except one, the process is discussed in  \cite{[1]}. For $C_n$ lattices $\mathcal{L}$ of rank $\geq 3$, we have seen that $A(\mathcal{L})=\sqcup_{i=1}^k (A_i\sqcup A_i^*)$. Therefore, we can use the same deletion process, deleting all pairs $a,a^*$  in every pair of atoms $A_i, A_i^*$ except one. 

In the following Theorems, we take independent sets as sets of atoms using the process above. In addition, using Lemmas \ref{Lemma9} and \ref{Lemma10} we can give a $C_n$ lattice structure to the restriction $[A,J]$ for every $A\in A(\mathcal{L})$. As the $*$ function on the restriction is induced by the original function on $\mathcal{L}$, we get the following relation between admissible independent sets:

\subsection{Lemma:}
\label{Lemma12}
Let $\omega$ be a linear order on the atoms of a $C_n$ lattice $\mathcal{L}$ of rank $d\geq 4$ with the first $d-1$ atoms forming a geometric independent set $\mathcal{I}$. Then for every $A_i\in A(\mathcal{L})$ there exists a linear order $\omega _i$ on the atoms of $\mathcal{L}_i=[A_i,\hat{1}]$ with the atoms that come first in the ordering being those that cover some $A_{j}$ where $j<i$ and the first $d-2$ atoms form a geometric independent set of $\mathcal{L}_i$.

\subsubsection*{Proof:}
Let $\omega$ be such an order and define the induced order $\omega_i$ on $A(\mathcal{L}_i)$ taking $B_j< B_k$ if $B_j$ covers an atom smaller than any atom covered by $B_k$ (excluding $A_i$) with respect to $\omega$. As the atoms of $\mathcal{L}_i$ partition $A(\mathcal{L})\setminus \left\{ A_{i},A_{i}^{*}\right\} $, $\omega _i$ is a well-defined linear order. We now observe two cases:
\begin{enumerate}
\item if $i\leq d-1$, then using Theorem \ref{Theorem8}, the first $d-2$ atoms of $\omega_i$ will be $A_1\vee A_i,...,A_{d-1}\vee A_i$. This is a geometric independent set of $\mathcal{L}_i$, and therefore $\omega_i$ is the required order.
\item If $i\geq d$, we again have two options. If $A_i^*\notin \vee \mathcal{I}$ then, using Theorem \ref{Theorem9} there exist $\mathcal{I}'\subsetneq \mathcal{I}$ such that $\mathcal{I}'\cup \left\{ A_{i}\right\}$ is a geometric independent set of size $d-1$. Otherwise, $A_i^*\in \vee \mathcal{I}$ and by Theorem \ref{Theorem9} there exist $\mathcal{I}'\subsetneq \mathcal{I}$ such that $\mathcal{I}'\cup \left\{ A_{i}^*\right\}$ is a geometric independent set of size $d-1$. applying Lemma \ref{Lemma8} we have $\mathcal{I}'\cup \left\{ A_{i}\right\}$ also a geometric independent set of size $d-1$.

In both cases we have any order $\omega_i'$ starting with $\left\{A_i\vee A_j\mid A_j\in \mathcal{I}'\right\}$ the required order.
\end{enumerate}
 
\subsection{Corollary:}
\label{corollary2}
If $\left\{ a_{1},...,a_{k}\right\}$ is a geometric independent set of a finite, graded, and atomic lattice $\mathcal{S}$, then $\left\{ a_{1}\vee a_{i},...,a_{k}\vee a_{i}\right\}\setminus \left\{ a_{i}\right\}  $ is a geometric independent set of $\left[a_{i},\hat{1}\right]$.

\subsection{Theorem:}
\label{Theorem14}
Every $C_n$ lattice $\mathcal{L}$ admits a recursive atom ordering.

\subsubsection*{Proof:}
First, notice that the conditions for a recursive atom ordering are met if the rank of $\mathcal{L}$ is $2$, since then all atoms are covered by $J$.
Another separated case is $rank(\mathcal{L})=3$, in this case taking any linear order $\omega$ on $A(\mathcal{L})$ with $A_1\neq {A_{2}}^*$ is a recursive atom ordering. The first condition is obvious, as all the elements of rank $2$ are covered by $J$. To see the second condition holds, observe that every pair of atoms which are not stars of each other is covered by a rank $2$. If a pair of atoms $A_i,A_j$ with $i\lneq j$ are stars of each other, then $\left\{i,j\right\}\neq \left\{1,2\right\}$. Therefore, $A_1\vee A_j$ covers $A_1,A,j$ or $A_2\vee A_j$ covers $A_2,A,j$.

We now prove the Theorem by induction on $t$, the number of atoms of $\mathcal{L}$. 
For the base case, we take $t=2$. This is a degenerate case
in which the atoms are also the maximal admissible elements. Assume
$A,B$ are atoms of $\mathcal{L}$, if there is some element $C\in\mathcal{L}$
that is not an atom we would have $C=A\vee B=A\cup B$ since $\mathcal{L}$
is atomic and since the set of atoms must partition $J$, we have
$C=J$.
\\

We now prove that any linear order on the atoms of a $C_n$ lattice $\mathcal{L}$ of rank $d$ with its first $d-1$ atoms forming a geometric independent set is a recursive atom ordering. We assume that the statement is true for every $C_n$ lattice $\mathcal{L}'$ with $\leq t-2$ atoms and prove it for a $C_n$ lattice $\mathcal{L}$ with $t$ atoms.

Let $A_{1},...,A_{t}$ be such an atom ordering. For the first condition, we just apply Lemma \ref{Lemma11} and Lemma \ref{Lemma12} to ensure that $\left[A_i,J\right]\subseteq\mathcal{L}$ has $\leq t-2$ atoms. Therefore, we get a recursive atom ordering using the induction hypothesis.

The proof of the second condition is the same as in the separate case $rank(\mathcal{L})=3$ as linear order $\omega_i$ starts with a geometric independent set which cannot contain an atom and its star..

\section{characterization by recursive atom ordering:}
In this section we prove a characterization of $C_n$
lattices by recursive atom ordering inspired by the work in \cite{[8]}, our goal is to give a minimal list of orders that are always recursive atom orderings equivalent to being a $C_n$ lattice. We start with a definition of strongly admissible sets characterizing non maximal geometric independent sets in $C_n$ lattices.

\subsection{Definition:}
\label{Definition17}
A set $A$ of atoms is said to be strongly admissible if either $A=\emptyset$, or there exists a strongly admissible subset $B\subsetneq A$ with $\left|B    \right|=\left|A\right|-1$ such that for the unique element $a\in A\setminus B$, we have $a^*\nless \vee B$.\\

We observe that every strongly admissible set is indeed admissible and that for any admissible element $X$ of a finite, atomic, and graded lattice $\mathcal{L}$ we have a spanning strongly admissible independent set. Next, we present a characterization of $C_n$ lattices similar to the characterization of geometric lattices given in Theorem \ref{Theorem6}.

\subsection{Lemma:}
\label{Lemma13}
Let $\mathcal{L}$ be a finite, atomic, and graded lattice. For every $a,b\in A(\mathcal{L})$ There exists a geometric independent set $a,b\in A$ such that $rank(a\vee b)=|I|$ .

\subsubsection*{Proof:}
Using Lemma \ref{Lemma4} we will construct a spanning geometric independent set for $a\vee b$, $B=\left\{a=a_1,a_2,...,a_{rank(a\vee b)}\right\}$ such that $n=rank(\vee_{i=1}^n a_i)$ for all $1 \leq n\leq rank(a\vee b)$ and $\vee B=a\vee b$. Let $1 \leq n\leq rank(a\vee b)$ be minimal so that $b\in \vee_{i=1}^n a_i$ and define $I=(B\setminus \left\{a_n\right\})\cup \left\{b\right\}$. We again have by Lemma \ref{Lemma4} that $I$ is independent and so our desired geometric independent set.

\subsection{Theorem:}
\label{Theorem15}
A finite, atomic, and graded lattice $\mathcal{L}$ is order isomorphic to a $C_n$ lattice if and only if there is a partition of the atoms of $\mathcal{L}$ to pairs ${1},{1}^{*},{2},{2}^{*},...,{k},{k}^{*}$ such that the non-maximal geometric independent sets are exactly the strongly admissible non maximal sets.

\subsubsection*{Proof:}
If $\mathcal{L}$ is a $C_n$ lattice, then equality is a consequence of Lemma \ref{Lemma7}.
Conversely, every atomic lattice $\left(\mathcal{L},\leq\right)$ is order isomorphic to a set of subsets,
ordered by inclusion by defining $x=\left\{ a\in A(\mathcal{L})\mid a\leq x\right\} $ for every $x\in \mathcal{L}$.
We prove the four conditions in Definition \ref{Definition3}:
\begin{enumerate}
\item As $\mathcal{L}$ is a finite lattice, we just have to denote the minimal set $\phi$ and the maximal set $J$. Notice that since $\mathcal{L}$ is atomic, we must have $J=\left[k\right]\sqcup\left[k\right]^{*}$.
\item Let $A\in \mathcal{L}$ be an inadmissible set, $\left\{ i,i^*\right\}\subset A$. Using Lemma \ref{Lemma13} we can extend $\left\{ i,i^*\right\}$ to a size $rank(A)$ independent set $I$ that spans $A$. As $I$ is inadmissible and geometric, we must have $rank(I)=rank(A)=rank(\mathcal{L})$ and so $A=J$. 
\item This property is always true for a finite atomic lattice of subsets. $A\wedge B\subseteq A\cap B$ and if there is an atom $i\in A\cap B$ not in the meet then $A\wedge B\subsetneq (A\wedge B)\vee \left\{ i\right\}\subset A,B$ which is a contradiction to the definition of $A\wedge B$, therefore $A\wedge B=A\cap B$.
\item Let $B,C\in \mathcal{L}$ be two elements that cover some $A\in \mathcal{L}$, since $B\cap C\in \mathcal{L}$ we must have $B\cap C=A$. Furthermore, if $A$ is not covered by $J$, then for $i\in A$ we have $I\vee i^*$ not geometric for every spanning geometric set $I$ of $A$ . Therefore, there is no element covering $A$ and containing $i^*$. It remains to prove that if $j\in J\setminus (A\sqcup A^*)$ then there exists $B\in \mathcal{L}$ covering $A$ with $j\in B$. Take a strongly admissible independent set $I$ that spans $A$, making $I\cup \left\{ i\right\}$ also a strongly admissible independent set. Therefore, $I\cup \left\{ i\right\}$ is geometric and we have $rank(I\cup \left\{ i\right\})= rank\left(A\right)+1$. The element of $\mathcal{L}$ spanned by $I\cup \left\{ i\right\}$ is the desired $B$.
\end{enumerate}

\subsection{Corollary:}
\label{corollary3}
A finite, atomic, and graded lattice $\mathcal{L}$ of rank $d$ is order isomorphic to a $C_n$ lattice if and only if there is a partition of the atoms of $\mathcal{L}$ to pairs ${1},{1}^{*},{2},{2}^{*},...,{k},{k}^{*}$ such that the size $d-1$ geometric independent sets are exactly the size $d-1$ strongly admissible sets.

\subsection{Lemma:}
\label{Lemma14}
Let $\mathcal{L}$ be a finite, atomic, and graded lattice, and let $I$ be a set of atoms of $\mathcal{L}$. If every ordering of the atoms of $\mathcal{L}$ starting with $I$ yields a recursive atom ordering, then $I$ is a geometric independent set.

\subsubsection*{Proof:}

We prove this by induction on the size of $I$. The base case $|I|=1$ is trivially geometric, since the rank of a single atom is 1.

Assume that every proper subset of $I=\{ a_1,a_2,\dots ,a_{|I|}\}$ is geometric. We aim to show that $I$ itself is geometric.
Let $\omega _1$ be a recursive atom ordering of $\mathcal{L}$ that begins with the atoms of $I$, i.e.,
\[
\omega _1=a_1\lneq a_2\lneq \cdots \lneq a_{|I|}\lneq \cdots 
\]
Now consider the induced recursive atom ordering $\omega _2$ on the interval $[a_{|I|},\hat {1}]$. This ordering begins with the set:
\[
I'=\{ a_1\vee a_{|I|},a_2\vee a_{|I|},\dots ,a_{|I|-1}\vee a_{|I|}\}
\]
According to the induction hypothesis, each pair $\{ a_i,a_{|I|}\}$  is geometric, so each join $a_i\vee a_{|I|}$ has rank 2. Moreover, these joins are distinct, since the join of three atoms is of rank $3$.
We repeat this process iteratively: for each $i=3,\dots ,|I|-1$, define $\omega _i$ as the recursive atom ordering on $[b_{|I|+1-i},\hat {1}]$, where $b_{|I|-i}$ is the $(|I|-i)$th atom in $\omega _{i-1}$.
The order $\omega _{|I|-1}$ begins with two elements of the form:
\[
\left\{ (\vee_{\underset{i\neq j,k}{1\leq i\leq\left|I\right|}}a_{i})\vee a_{j},(\vee_{\underset{i\neq j,l}{1\leq i\leq\left|I\right|}}a_{i})\vee a_{j}\right\} 
\]
for distinct $j,k,l\in I$.

Since $\omega _{|I|-1}$ is a recursive atom ordering, the join of these two elements must have rank $2$ in the interval $[\vee_{\underset{i\neq l,k}{1\leq i\leq\left|I\right|}}a_{i},\hat{1}]$. Thus, the full join satisfies:
\[
rank(\vee_{1\leq i\leq\left|I\right|}a_i)=rank(((\vee_{\underset{i\neq j,k}{1\leq i\leq\left|I\right|}}a_{i})\vee a_{j})\vee((\vee_{\underset{i\neq j,l}{1\leq i\leq\left|I\right|}}a_{i})\vee a_{j}))
\]
\[
=rank((\vee_{\underset{i\neq j,k}{1\leq i\leq\left|I\right|}}a_{i})\vee a_{j})+1=|I|
\]
Hence, I is a geometric independent set.\\

We are now ready to prove the main Theorem of this section.

\subsection{Theorem:}
\label{Theorem16}
Let $\mathcal{L}$ be a finite, atomic, and graded lattice of rank $3\leq d\lneq 2k=\left|A\left(\mathcal{L}\right)\right|$. then $\mathcal{L}$ is order isomorphic to a $C_n$ lattice if and only if the following condition holds:

There exists a partition of the atom set $A(\mathcal{L})$ into $k$ disjoint pairs:
\[
\left\{{1},{1}^{*}\right\},\left\{{2},{2}^{*}\right\},...,\left\{{k},{k}^{*}\right\}
\]
such that every linear ordering
\[
a_1\lneq ...\lneq a_{2k}
\]
of the atoms beginning with a subset of size $d-1$ that is either strongly admissible or geometric independent, constitutes a recursive atom ordering if and only if the condition
\[
a_i\neq {a_{i+1}}^* \textbf{ for all  } i\in [d-2]
\]
is satisfied.

\subsubsection*{Proof:}

If $\mathcal{L}$ is a $C_n$ lattice, every strongly admissible independent set of size $\leq d-1$ is geometric, and we have proven in Theorem \ref{Theorem14} that such a recursive atom ordering exists. 

To prove the converse, we invoke Corollary \ref{corollary3}. Throughout the proof, we assume $d\geq 4$, since the case $d=3$ is straightforward.
When $d=3$, condition (1) of Definition \ref{Definition16} holds trivially, and condition (2) immediately implies that admissible pairs of atoms have rank $2$, while inadmissible pairs have rank $3$.

Lemma \ref{Lemma14} states that if every ordering of the atoms of $\mathcal{L}$ starting with a set $I$ yields a recursive atom ordering, then I is a geometric independent set. Therefore, every strongly admissible set is a geometric independent set.

We now aim to prove that $rank(i\vee i^*)=d$ for every $i\in [k]$. This result implies that any geometric independent set of size less than $d$ must be strongly admissible.
Suppose, for contradiction, that this is not the case. Then, by Lemma \ref{Lemma13}, there exists a geometric independent set $A\subseteq A(\mathcal{L})$ of size $d-1$ that includes both $i$ and $i^*$. Under this assumption, we obtain a recursive atom ordering $\omega$  of the form:
\[
i\prec j\prec i^*\prec \cdots 
\]
for some $j\in A$. It is straightforward to verify that the modified ordering $\omega _{i,i^*}$, obtained by swapping the positions of $i$ and $i^*$, also satisfies the conditions of a recursive atom ordering.

We now verify that condition (1) of Definition \ref{Definition16} holds for the ordering $\omega _{i,j}$, which is obtained from a recursive atom ordering $\omega$  by swapping the positions of $i$ and $j$.
Since $\omega$  is a recursive atom ordering, it suffices to show that the interval $[i,\hat {1}]$ admits a recursive atom ordering that begins with $i\vee j$.
\begin{itemize}
    \item {Case 1: $d\geq 5$} 
    
        In this case, we can construct a recursive atom ordering of $\mathcal{L}$ that begins with the sequence $j\prec i\prec t\prec i^*$ for some atom $t\in I$. This ordering induces a recursive atom ordering on $[i,\hat {1}]$ that starts with $i\vee j$, as required.
    \item {Case 2: $d=4$}
    
        Here, it suffices to find an atom $t$ such that $\{ i,j,t\}$  forms a strongly admissible independent set. Such a set exists unless $i\vee j=\bigvee T$ for some transversal $T$. However, by symmetry, this would also imply $i^*\vee j=\bigvee T$, leading to the conclusion that the atom set is $A(\mathcal{L})=\{ i,i^*,j,j^*\} $ which contradicts the assumption that $rank(\mathcal{L})<|A(\mathcal{L})|$.
\end{itemize} 

Since condition (1) of Definition \ref{Definition16} is satisfied for the ordering $\omega _{i,j}$, yet $\omega _{i,j}$ fails to be a recursive atom ordering, it follows that condition (2) must be violated. Consequently, we deduce that:
\[
rank(i\vee i^*)\gneq 2 \textbf{ and } j\vee i^*\nleq i\vee i^*
\]

This contradicts condition (1) of Definition \ref{Definition16} as applied to the ordering $\omega$.\\
\\
\textbf{Acknowledgements.} The author is supported by
Horizon Europe ERC Grant number: 101045750 / Project acronym: HodgeGeoComb.

\end{document}